\documentclass[final]{amsart}

\usepackage{amsthm,amsfonts,amssymb,amsmath,latexsym, amscd, mathrsfs, oldgerm,float}
\usepackage{exscale, cite, epsfig, graphics}
\usepackage{color}
\usepackage[dvipsnames]{xcolor}
\usepackage{tikz}

\usepackage{accents}

\usepackage[notref, notcite]{showkeys}
\usepackage[colorlinks]{hyperref}

\renewcommand{\leq}{\leqslant}

\newtheorem*{main-theorem}{Main Theorem}
\newtheorem*{remark*}{Remark}
\numberwithin{equation}{section}

\title[Shallow water models with constant vorticity]
{Shallow water models with constant vorticity}

\author[Hur]{Vera~Mikyoung~Hur}
\address{Department of Mathematics, University of Illinois at Urbana-Champaign, Urbana, IL 61801 USA}
\email{verahur@math.uiuc.edu}

\date{\today}


\begin{document}

\maketitle

\begin{abstract}
We modify the nonlinear shallow water equations, the Korteweg-de Vries equation, and the Whitham equation, to permit constant vorticity, and examine wave breaking, or the lack thereof. By wave breaking, we mean that the solution remains bounded but its slope becomes unbounded in finite time. We show that a solution of the vorticity-modified shallow water equations breaks if it carries an increase of elevation; the breaking time decreases to zero as the size of vorticity increases. We propose a full-dispersion shallow water model, which combines the dispersion relation of water waves and the nonlinear shallow water equations in the constant vorticity setting, and which extends the Whitham equation to permit bidirectional propagation. We show that its small amplitude and periodic traveling wave is unstable to long wavelength perturbations if the wave number is greater than a critical value, and stable otherwise, similarly to the Benjamin-Feir instability in the irrotational setting; the critical wave number grows unboundedly large with the size of vorticity. The result agrees with that from a multiple scale expansion of the physical problem. We show that vorticity considerably alters the modulational stability and instability in the presence of the effects of surface tension.
\end{abstract}

\section{Introduction}\label{sec:intro}

Much of the theory of surface water waves makes the assumption that the flow is irrotational. This is well justified in some circumstances. Moreover, in the absence of the initial vorticity, boundaries, or currents, water waves will have no vorticity at all future times. On the other hand, rotational effects are significant in many circumstances. For instance, in any region where wind is blowing, there is a surface drift of the water, and wave parameters, such as maximum wave height, are sensitive to the velocity at a wind-drift boundary layer. Moreover, currents cause shear at the bed of the sea or a river; see \cite{SilvaPeregrine1988}, for instance.

In 1802 Gerstner~\cite{Gerstner} (see also \cite{Constantin2001}) produced a remarkable example of periodic traveling waves in water of infinite depth for a certain nonzero vorticity. No solution formula of the kind is known in the irrotational setting. Perhaps more striking are internal stagnation points and critical layers; see \cite{RMN2017}, for instance, and references therein. By the way, the flow beneath a periodic traveling wave is necessarily supercritical in the irrotational setting. Moreover, vorticity influences the stability and instability of a shear flow, with or without free surface; see \cite{HL2008}, for instance. For more about wave current interactions, we encourage the interested reader to an excellent textbook \cite{Constantin-book}.

Here we restrict the attention to constant vorticity. This is interesting because of analytical tractability. As a matter of fact, the governing equations of the water wave problem may be written in terms of quantities at the fluid surface, similarly to the irrotational setting. Moreover, constant vorticity is representative of a wide range of physical scenarios. When waves are short compared with the vorticity length scale, the vorticity at the surface layer is dominant in the wave dynamics. Moreover, when waves are long compared with the fluid depth, the mean vorticity is more important than its specific distribution; see \cite{SilvaPeregrine1988}, for instance. Examples include tidal currents --- alternating horizontal movements of water associated with the rise and fall of the tide --- where positive and negative vorticities are appropriate for ebbs and floods, respectively; see \cite{Constantin-book}, for instance. 

It is a matter of experience that the smooth surface of a water wave may transform into rough fronts of spray and bubbles; see \cite{Peregrine1983}, for instance. A direct study of ``breaking" is difficult because the governing equations of the water wave problem are very complicated. Vorticity makes it even more difficult. One may resort to simpler approximate models to gain insights. 

Here we explore shallow water models in the presence of constant vorticity. By shallow water, we mean that waves are long compared with the fluid depth. Specifically, we modify the nonlinear shallow water equations, the Korteweg-de Vries equation, and the Whitham equation, to permit constant vorticity; see \eqref{E:shallow}, \eqref{E:KdV}, and \eqref{E:Whitham}. One may follow along the same line as the argument in \cite{Lannes}, for instance, in the irrotational setting, to rigorously justify them. After submitting the manuscript, the author has learned that Kharif and Abid~\cite{KA2017} derived them independently. For general vorticities, Freeman and Johnson~\cite{FJ1970} derived the Korteweg-de Vries equation. For constant vorticities, Johnson~\cite{Johnson2012} developed the Boussinesq and Camassa-Holm equations. Moreover, for general vorticities, Castro and Lannes~\cite{CL2014} recently developed equations of Green-Naghdi type.

For any constant vorticity, we show that a solution of \eqref{E:shallow} carrying an increase of elevation breaks, whereas no solution of \eqref{E:KdV} breaks, similarly to the irrotational setting (see \cite{Whitham}, for instance). By wave breaking, we mean that the solution remains bounded but its slope becomes unbounded in finite time. The breaking time decreases to zero as the size of vorticity increases. In particular, for any $t>0$, a nontrivial solution of \eqref{E:shallow} breaks at the time $t$, provided that the size of vorticity is sufficiently large, depending on $t$. For any constant vorticity, we extend the proof in \cite{Hur-breaking} and show that a solution of \eqref{E:Whitham} breaks, provided that the initial datum is sufficiently steep, similarly to the irrotational setting (see \cite{Hur-breaking}).

Moreover, we propose ``full-dispersion shallow water equations" with constant vorticity, which combine \eqref{E:shallow} and the dispersion relation in the linear theory of water waves, and which extend \eqref{E:Whitham} to permit bidirectional propagation. For any constant vorticity, we show that a small amplitude and periodic traveling wave of \eqref{E:FDSW} is unstable to long wavelength perturbations, provided that the wave number is greater than a critical value, and stable otherwise, similarly to the irrotational setting (see \cite{HP2}). 
The critical wave number grows unboundedly large with the size of vorticity. In particular, for any $k>0$, a small amplitude and $2\pi/k$ periodic traveling wave of \eqref{E:FDSW} is modulationally stable, provided that the size of vorticity is sufficiently large, depending on $k$. 
The result qualitatively agrees with that in \cite{TKM2012}, for instance, from a multiple scale expansion of the physical problem (see Figure~\ref{fig:T=0'}), and it improves that in \cite{HJ3} based on \eqref{E:Whitham}. Moreover, constant vorticity considerably alters the modulational stability and instability in \eqref{E:FDSW} in the presence of the effects of surface tension (see Figure~\ref{fig:w=0}, Figure~\ref{fig:w=-3}, and Figure~\ref{fig:w=3}). It is interesting to compare the result with that from a multiple scale expansion of the physical problem.

Of course, at the most prominent stage of breaking, an element of the fluid surface becomes vertical; a portion of the surface overturns, projects forward, and forms a jet of water; see \cite{Peregrine1983}, for instance. Overturning Stokes waves with constant vorticity will be numerically studied in \cite{DH1}. By the way, the profile of a periodic traveling wave is necessarily the graph of a single valued function in the irrotational setting.

\section{The water wave problem with constant vorticity}\label{sec:WW}

The water wave problem, in the simplest form, concerns the wave motion at the surface of an incompressible and inviscid fluid, lying below a body of air, and under the influence of gravity. Although an incompressible fluid such as water may have variable density, we assume for simplicity that the density $=1$. Suppose for definiteness that in Cartesian coordinates, waves propagate in the $x$ direction and gravity acts in the negative $y$ direction. 
Suppose that the fluid at rest occupies the region, bounded above by the free surface $y=0$ and below by the rigid bottom $y=-h_0$ for some constant $h_0>0$. 

Let $y=\eta(x;t)$ represent the fluid surface displacement from $y=0$ at the point $x$ and the time $t$; $h_0+\eta>0$ is physically realistic. Let 
\[
(u(x,y;t)-\omega y, v(x,y;t))
\] 
denote the velocity of the fluid at the point $(x,y)$ and the time $t$ for some constant $\omega\in\mathbb{R}$, and $p(x,y;t)$ the pressure. In the bulk of the fluid, they satisfy the Euler equations for an incompressible fluid:
\begin{equation}\label{E:euler}
\begin{aligned}
&u_t+(u-\omega y)u_x+v(u_y-\omega)=-p_x, \\
&v_t+(u-\omega y)v_x+vv_y=-p_y-g, \\
&u_x+v_y=0,
\end{aligned}
\end{equation}
where $g>0$ is the constant due to gravitational acceleration. We assume that 
\begin{equation}\label{E:vorticity}
v_x-u_y=0
\end{equation}
for $-h_0<y<\eta(x;t)$. In other words, $\omega$ is the constant vorticity throughout the fluid region. By the way, if vorticity is constant at $t=0$, then it remains so at all future times, as long as the fluid region is two dimensional and simply connected. 

The kinematic and dynamic conditions at the fluid surface:
\begin{equation}\label{E:surface}
v=\eta_t+(u-\omega \eta)\eta_x\quad\text{and}\quad p=p_{atm}
\end{equation}
state, respectively, that fluid particles do not invade the air, nor vice versa, and that the pressure at the fluid surface equals the constant atmospheric pressure $p_{atm}$. This neglects the motion of the air and the effects of surface tension (see \eqref{E:T-surface}, for instance). The kinematic condition at the bottom:
\begin{equation}\label{E:bottom}
v=0
\end{equation}
states that the fluid particles at the bottom remain so at all times. 

Thanks to the last equation of \eqref{E:euler} and \eqref{E:vorticity}, we may define a velocity potential for the irrotational perturbation as
\begin{equation}\label{def:phi}
\nabla \phi=(u,v).
\end{equation}
Note that $\Delta\phi=0$ inside the fluid. Let $\psi$ be a harmonic conjugate of $\phi$. We may rewrite \eqref{E:surface} in terms of $\phi$ (and $\psi$) as
\[
\begin{aligned}
&\eta_t+(\phi_x-\omega\eta)\eta_x=\phi_y, \\
&\phi_t+\frac12|\nabla\phi|^2-\omega\eta\phi_x+\omega\psi+g\eta=0,
\end{aligned}
\]
say. We may then rewrite \eqref{E:euler}-\eqref{E:bottom} in terms of $\eta(x;t)$ and $\phi(x,\eta(x;t);t)$ at the fluid surface; see \cite{DH1}, for instance, for details.

\subsection*{Trivial solution}

For any $h_0>0$ and $\omega\in\mathbb{R}$, note that
\begin{equation}\label{E:trivial}
\eta(x;t)=0, \qquad (u,v)(x,y;t)=(0,0),\quad\text{and}\quad p(x,y;t)=p_{atm}-gy,
\end{equation}
for $t,x\in\mathbb{R}$ and $-h_0\leq y\leq 0$, solve \eqref{E:euler}-\eqref{E:bottom}. Physically, they make a shear flow in a channel of depth $h_0$, for which the fluid surface and the velocity are horizontal, $\omega$ is the constant vorticity, and the pressure is hydrostatic. Here we restrict the attention to waves propagating in a shear flow of the kind.

\subsection*{Stokes waves}

It is a matter of experience that waves commonly seen in the ocean or a lake are approximately periodic and propagate over a long distance practically at a constant speed without change of form. Stokes in his 1847 memoir (see also \cite{Stokes}) made many contributions about waves of the kind, for instance, observing that crests would be sharper and throughs flatter as the amplitude increases, and a wave of greatest possible height would exhibit stagnation with a $120^\circ$ corner. By stagnation, we mean that the fluid particle velocity coincides with the wave speed. This is argued for incipient breaking; see \cite{Peregrine1983}, for instance. In the irrotational setting, the rigorous existence theory of Stokes waves is nearly complete. We encourage the interested reader to excellent surveys~\cite{Toland1996, BT-book}. On the other hand, for general vorticities, it was not until when Constantin and Strauss~\cite{CS2004} established the existence from no wave up to, but not including, a ``limiting" wave with stagnation. 

For constant vorticities, Constantin, Strauss, and Varvaruca~\cite{CSV} recently made strong use of that $\phi$ in \eqref{def:phi} is harmonic inside the fluid, whereby they followed along the same line as an argument in the irrotational setting (see \cite{BT-book}, for instance) for the global bifurcation of periodic traveling waves, permitting internal stagnation, critical layers, and overturning surface profiles. As a matter of fact, numerical computations (see \cite{SilvaPeregrine1988, RMN2017}, for instance) confirm critical layers and overturning profiles. By the way, in the irrotational setting, the flow beneath a Stokes wave is necessarily supercritical and the profile is the graph of a single valued function. Moreover, Constantin, Strauss, and Varvaruca conjectured that either a limiting wave would exhibit stagnation with a $120^\circ$ corner, or its surface profile would overturn and intersect itself at the trough line. This will be numerically studied in \cite{DH1}. 

\section{The linear problem}\label{sec:linear}

We linearize \eqref{E:euler}-\eqref{E:bottom} about \eqref{E:trivial} to arrive at
\[
\left\{\begin{split}
&u_t-\omega yu_x-\omega v=-p_x, \\
&v_t-\omega yv_x=-p_y, \\
&u_x+v_y=0\quad\text{and}\quad v_x-u_y=0
\end{split}\right.
\]
in $-h_0<y<0$, and
\begin{align*}
&v=\eta_t\quad\text{and}\quad p=g\eta &&\text{at}\quad y=0, \\
&v=0 &&\text{at}\quad y=-h_0.
\end{align*}
Seeking a solution, $2\pi/k$ periodic and traveling at the speed $c$ in the $x$ direction, it is reasonable to assume that
\[
\eta(x;t)=\cos(k(x-ct)).
\]
A straightforward calculation reveals that 
\begin{align*}
&u(x,y;t)=\cos(k(x-ct))\frac{ck\cosh(k(h_0+y))}{\sinh(kh_0)},\\
&v(x,y;t)=\sin(k(x-ct))\frac{ck\sinh(k(h_0+y))}{\sinh(kh_0)}, \\
\intertext{and}
&p(x,y;t)=\cos(k(x-ct))\frac{(c+\omega y)ck\cosh(k(h_0+y))-\omega c\sinh(k(h_0+y))}{\sinh(kh_0)},
\end{align*}
where $c$ satisfies
\begin{equation}\label{E:dispersion}
\Big(c-\frac{\omega}{2}\frac{\tanh(kh_0)}{k}\Big)^2=
\frac{g\tanh(kh_0)}{k}+\frac{\omega^2}{4}\frac{\tanh^2(kh_0)}{k^2}
\end{equation}
--- namely, the dispersion relation. Physically, $c$ means the speed in the linear theory of a $2\pi/k$ periodic wave traveling at the surface of a shear flow with the constant vorticity $\omega$ in a channel of depth $h_0$. Note that a positive vorticity increases the phase speed, while a negative vorticity decreases it. In the irrotational setting, \eqref{E:dispersion} becomes the well-known dispersion relation 
\begin{equation}\label{E:dispersion0}
c^2=\frac{g\tanh(kh_0)}{k};
\end{equation}
compare \cite[(13.25)]{Whitham}, for instance. 

For general vorticities, the dispersion relation may not be determined explicitly. Rather, it is based on a certain Sturm-Liouville problem; see \cite{CS2004, HL2008, Kar12}, for instance, for details. 

A critical layer occurs if $c=-\omega y$ for some $y$ in the range $(-h_0,0)$. Note from \eqref{E:dispersion} that it does when
\[
\omega^2>\frac{g\tanh(kh_0)}{kh_0^2-h_0\tanh(kh_0)}.
\]

In the limit as $kh_0\to0$, one may approximate \eqref{E:dispersion} by
\begin{equation}\label{E:dispersion1}
\Big(c-\frac12\omega h_0\Big)^2\sim gh_0+\frac14\omega^2h_0^2.
\end{equation}
Note that the limiting phase speed is independent of $k$. In the irrotational setting, \eqref{E:dispersion1} corresponds to the critical Froude number. 

Moreover, in the zero gravity limit as $g\to0$, one may approximate \eqref{E:dispersion} by 
\begin{equation}\label{E:g=0}
c\sim\omega\frac{\tanh(kh_0)}{k}
\end{equation}
for $\omega>0$ and $c\sim0$ for $\omega<0$.

\subsection*{The effects of surface tension}\label{sec:surface tension}

Surface tension acts like a stretched membrane at the fluid surface. We continue to assume that the fluid is incompressible inviscid, the vorticity is constant, and we use the notation of Section~\ref{sec:WW}. The main difference is that the pressure at the fluid surface no longer equals the atmospheric pressure. Instead, the jump of the pressure across the fluid surface is proportional to the curvature of the surface. Therefore,
\begin{equation}\label{E:T-surface}
p=p_{atm}-T\Big(\frac{1}{\sqrt{1+\eta_x^2}}\Big)_x \qquad\text{at}\quad y=\eta(x;t)
\end{equation}
replaces the latter equation of \eqref{E:surface}, where $T>0$ is the coefficient of surface tension. 
We repeat the preceding argument and find the dispersion relation 
\begin{equation}\label{E:T-dispersion}
\Big(c-\frac{\omega}{2}\frac{\tanh(kh_0)}{k}\Big)^2
=(g+Tk^2)\frac{\tanh(kh_0)}{k}+\frac{\omega^2}{4}\frac{\tanh^2(kh_0)}{k^2}.
\end{equation}
In the irrotational setting, \eqref{E:T-dispersion} becomes
\begin{equation}\label{E:T-dispersion0}
c^2=(g+Tk^2)\frac{g\tanh(kh_0)}{k};
\end{equation}
compare \cite[(13.52)]{Whitham}, for instance.

At air-sea interfaces, $g\approx 9.81 m/s^2$ and $T\approx 7.3\times 10^{-3} N/m$. Therefore, for any $\omega\in\mathbb{R}$, the effects of surface tension become negligible if wavelengths are much greater than $2\pi\sqrt{T/g}\approx 1.7cm$, for instance, in a plunging jet and the entrainment of air. On the other hand, we learn in Section~\ref{sec:MI} that the effects of surface tension considerably alter the modulational stability and instability of a periodic traveling wave, with or without vorticity. Moreover, gravity capillary waves, right down to the amplitude of $2mm$ or so, may develop strongly turbulent regions; see \cite{Peregrine1983}, for instance.

\section{The nonlinear shallow water equations}\label{sec:shallow}

We follow along the same line as the argument in \cite[Section~13.10]{Whitham} to approximately describe finite amplitude and shallow water waves of \eqref{E:euler}-\eqref{E:bottom} in the presence of constant vorticity. By shallow water, we mean that waves are long compared with the fluid depth. This will be made precise in \eqref{def:shallow water}.

We begin by approximating the second equation of \eqref{E:euler} by
\[
0=-p_y-g. 
\]
An integration leads to 
\begin{equation}\label{E:shallow-p}
p=p_{atm}+g(\eta-y).
\end{equation}
In other words, the pressure is hydrostatic at leading order. The first equation of \eqref{E:euler} becomes
\[
u_t+(u-\omega y)u_x+v(u_y-\omega)=-g\eta_x.
\]
Since the right hand side is independent of $y$, the rate of change of $u$ along each particle trajectory is independent of $y$. In particular, if $u$ is independent of $y$ throughout the fluid region at $t=0$, it remains so at all future times. Assume that $u$ does not depend on $y$ at $t=0$ at leading order. We evaluate the left hand side at the fluid bottom $y=-h_0$ and use \eqref{E:bottom}, to arrive at
\[
u_t+(u+\omega h_0)u_x=-g\eta_x,
\]
where $u=u(x,-h_0;t)$.

Although we reject $v$, and $v_t$, $v_x$ in the second equation of \eqref{E:euler} with respect to the other terms, we may not reject $v_y$ in the last equation of \eqref{E:euler}. Instead, we calculate
\begin{align*}
0=\int^{\eta(x;t)}_{-h_0}(u_x+v_y)(x,y;t)&~dy \\
=\Big(\int^{\eta(x;t)}_{-h_0} u(x,y;t)~dy\Big)_x&-u(x,\eta(x;t);t)\eta_x(x;t) \\
&+v(x,\eta(x;t);t)-v(x,-h_0;t) \\
=\Big(\int^{\eta(x;t)}_{-h_0} u(x,y;t)~dy\Big)_x&+\eta_t-\omega\eta\eta_x.
\end{align*}
Here the first equality uses the last equation of \eqref{E:euler}, the second equality uses the chain rule, and the last equality uses the former equation of \eqref{E:surface} and \eqref{E:bottom}.
Recall that $u$ does not depend on $y$ at leading order throughout the fluid region at all times. We evaluate the right hand side at the fluid bottom $y=-h_0$, to arrive at 
\[
\eta_t-\omega\eta\eta_x+(u(h_0+\eta))_x=0,
\]
where $u=u(x,-h_0;t)$. Together, 
\begin{equation}\label{E:shallow}
\begin{aligned}
&h_t+\omega h_0h_x+(uh)_x-\omega hh_x=0,\\
&u_t+\omega h_0u_x+gh_x+uu_x=0
\end{aligned}
\end{equation}
make {\em the vorticity-modified shallow water equations}, where
\begin{equation}\label{def:h}
h=h_0+\eta
\end{equation}
is the total depth of the fluid; $h>0$ is physically realistic. In the irrotational setting, \eqref{E:shallow} becomes 
\begin{equation}\label{E:shallow0}
\begin{aligned}
&h_t+(uh)_x=0,\\
&u_t+gh_x+uu_x=0,
\end{aligned}
\end{equation}
and agrees with \cite[(13.79)]{Whitham}, for instance.

Let's make an order of magnitude calculation in the preceding approximation. Note that the error for $p$ in \eqref{E:shallow-p} is of order $h_0v_t$. Since $v\approx -h_0u_x$ by the last equation of \eqref{E:euler}, the relative error in the latter equation of \eqref{E:shallow} is of order
\[
\frac{-p_x}{u_t}\approx \frac{h_0^2u_{xxt}}{u_t}\approx \frac{h_0^2}{\ell^2},
\]
where $\ell$ is the characteristic length scale in the $x$ direction. Therefore, one may regard \eqref{E:shallow} as an approximate model of \eqref{E:euler}-\eqref{E:bottom} for relatively shallow water or, equivalently, relatively long waves, so that 
\begin{equation}\label{def:shallow water}
(h_0/\ell)^2\ll 1.
\end{equation}
One may modify the argument in \cite[Section~6.2]{Lannes}, for instance, in the irrotational setting, to rigorously justify \eqref{E:shallow}. But we do not include the details here.

Note that the phase speed of the linear part of \eqref{E:shallow} satisfies
\[
\Big(c-\frac12\omega h_0\Big)^2=gh_0+\frac14\omega^2h_0^2.
\]
This agrees with \eqref{E:dispersion1}. It is independent of the wave number. In other words, the dispersion effects drop out in \eqref{E:shallow}. Of course, the dispersion effects are important in many circumstances, in particular, in coastal oceanography. In the following sections, we discuss how to include them.

\subsection*{Solutions of \eqref{E:shallow}}

We employ theory of hyperbolic conservation laws, solve \eqref{E:shallow} analytically, and discuss wave breaking. 

A straightforward calculation reveals that the characteristic velocities of \eqref{E:shallow} are 
\[
u-\omega \eta+\frac12\omega h\pm\sqrt{gh+\frac14\omega^2h^2}, 
\]
and the Riemann invariants of \eqref{E:shallow} are 
\begin{equation}\label{E:RI}
u-\frac12\omega h
\pm\bigg(\sqrt{gh+\frac14\omega^2h^2}
+\frac{g}{\omega}\log\bigg(2g+\omega^2h+2\omega\sqrt{gh+\frac14\omega^2h^2}\bigg)\bigg).
\end{equation}
In the irrotational limit as $\omega\to0$, they tend to $u\pm\sqrt{gh}$ and $u\pm2\sqrt{gh}$, respectively, and agree with the characteristic velocities and the Riemann invariants of \eqref{E:shallow0}; see \cite[Section~13.10]{Whitham}, for instance. Therefore, it follows from the method of characteristics (see \cite[Section~5.3]{Whitham}, for instance) that a simple wave solution of \eqref{E:shallow} propagating to the right into a shear flow with the constant vorticity $\omega$ in a channel of depth $h_0$ is
\begin{equation}\label{E:shallow solution}
\begin{aligned}
h(x,t)=&f(s), \\
u(x,t)=&\frac12\omega f(s)+\sqrt{gf(s)+\frac14\omega^2f^2(s)} \\
&+\frac{g}{\omega}\log\bigg(2g+\omega^2f(s)+2\omega\sqrt{gf(s)+\frac14\omega^2f^2(s)}\bigg) \\
&-\frac12\omega h_0
-\sqrt{gh_0+\frac14\omega^2h_0^2}
-\frac{g}{\omega}\log\bigg(2g+\omega^2h_0+2\omega\sqrt{gh_0+\frac14\omega^2h_0^2}\bigg),
\end{aligned}
\end{equation}
where 
\begin{equation}\label{E:char}
\begin{aligned}
x(t;s)=s+t\bigg(&\frac12\omega h_0+2\sqrt{gf(s)+\frac14\omega^2f^2(s)} \\
&+\frac{g}{\omega}\log\bigg(2g+\omega^2f(s)+2\omega\sqrt{gf(s)+\frac14\omega^2f^2(s)}\bigg) \\
&-\sqrt{gh_0+\frac14\omega^2h_0^2}
-\frac{g}{\omega}\log\bigg(2g+\omega^2h_0+2\omega\sqrt{gh_0+\frac14\omega^2h_0^2}\bigg)\bigg).
\end{aligned}
\end{equation}
In the irrotational limit as $\omega\to 0$, \eqref{E:shallow solution} and \eqref{E:char} tend to 
\begin{align*}
&h(x,t)=f(s),\qquad u(x,t)=2\sqrt{gh_0}-2\sqrt{gf(s)},\\
&x=s+(3\sqrt{gf(s)}-2\sqrt{gh_0})t,
\end{align*}
respectively, and agree with \cite[(13.80)]{Whitham}, for instance.

Note from \eqref{E:shallow solution}-\eqref{E:char} that for any $\omega\in\mathbb{R}$, a smooth solution of \eqref{E:shallow} carrying an increase of elevation and propagating into the undisturbed fluid depth breaks. By wave breaking, we mean that the solution remains bounded but its slope becomes unbounded in finite time. As a matter of fact, a straightforward calculation reveals that
\[
\frac{\partial x}{\partial s}(t;s)=1+\frac14\frac{3g+\omega^2f(s)}{\sqrt{gf(s)+\frac14\omega^2f^2(s)}}f'(s)t.
\]
Since $f(s)>0$ everywhere in $\mathbb{R}$, $\frac{\partial x}{\partial s}(t;s)=0$ for some $t>0$; the characteristics cross and the solution breaks. This occurs for the first time at 
\begin{equation}\label{def:T}
t_*=-\inf_{s\in\mathbb{R}}\frac{4}{f'(s)}\frac{\sqrt{gf(s)+\frac14\omega^2f^2(s)}}{3g+\omega^2f(s)}.
\end{equation}

Constant vorticity does not qualitatively change wave breaking in the shallow water theory. But the breaking time, defined in \eqref{def:T}, decreases to zero as the size of vorticity increases. In particular, for any $t_*>0$, the solution of \eqref{E:shallow} breaks at the time $t_*$, provided that $|\omega|$ is sufficiently large, depending on $t_*$.

\section{The Korteweg-de Vries equation}\label{sec:KdV}

We follow along the same line as the argument in \cite[Section~13.11]{Whitham}, combine some dispersion effects and \eqref{E:shallow}, to approximately describe small amplitude and shallow water waves of \eqref{E:euler}-\eqref{E:bottom} in the presence of constant vorticity. This will be made precise in \eqref{def:KdV regime}. 
 
We begin by linear waves propagating to the right at the surface of a shear flow with the constant vorticity $\omega$ in a channel of depth $h_0$. For shallow water waves satisfying \eqref{def:shallow water} or, equivalently, $kh_0\ll 1$, where $k$ is the wave number, we expand \eqref{E:dispersion} up to terms of order $(kh_0)^2$ to find
\begin{align}
c\sim&
\frac12\omega h_0+\sqrt{gh_0+\frac14\omega^2h_0^2}
-\frac16\Bigg(\omega h_0+\frac{gh_0+\frac12\omega^2h_0^2}{\sqrt{gh_0+\frac14\omega^2h_0^2}}\Bigg)(kh_0)^2 \notag\\
=:&c_0-c_2 k^2.\label{E:dispersion2}
\end{align}
Note that $c_0$ agrees with \eqref{E:dispersion1} and $c_2>0$ for any $\omega\in\mathbb{R}$. In the irrotational setting, 
\[
c_0=\sqrt{gh_0}\quad\text{and}\quad c_2=c_0h_0^2,
\]
and they agree with \cite[(13.94)]{Whitham}, for instance. A simplest partial differential equation whose dispersion relation is \eqref{E:dispersion2} would be
\begin{equation}\label{E:KdV-linear}
\eta_t+c_0\eta_x+c_2\eta_{xxx}=0.
\end{equation}

On the other hand, recall from the previous section that in the nonlinear shallow water theory, waves propagating to the right at the surface of a shear flow with the constant vorticity $\omega$ in a channel of depth $h_0$ satisfy the Riemann invariants \eqref{E:RI}. We substitute it in the former equation of \eqref{E:shallow} (or the latter equation), recall \eqref{def:h}, and make an explicit calculation, to arrive at
\[
\begin{aligned}
\eta_t-\omega\eta\eta_x+\bigg(\omega h&+2\sqrt{gh+\frac14\omega^2h^2}-\sqrt{gh_0+\frac14\omega^2h_0^2} \\
&+\frac{g}{w}\log\Big(2g+\omega^2h+2\omega\sqrt{gh+\frac14\omega^2h^2}\bigg) \\
&-\frac{g}{w}\log\Big(2g+\omega^2h_0+2\omega\sqrt{gh_0+\frac14\omega^2h_0^2}\bigg)\bigg)\eta_x=0.
\end{aligned}
\]
In the irrotational limit as $\omega\to0$, this tends to 
\[
\eta_t+(3\sqrt{gh}-2\sqrt{gh_0})\eta_x=0,
\]
and agrees with \cite[(13.97)]{Whitham}, for instance. Moreover, for small amplitude waves satisfying
\[
a/h_0\ll 1,\qquad\text{where $a$ is a typical amplitude,}
\]
we expand the nonlinearity up to terms of order $a/h_0$ to find
\begin{equation}\label{E:KdV-nonlinear}
\eta_t+\Bigg(\frac12\omega h_0+\sqrt{gh_0+\frac14\omega^2h_0^2}
+\frac12\frac{3g+\omega^2h_0}{\sqrt{gh_0+\frac14\omega^2h_0^2}}\eta\Bigg)\eta_x=0.
\end{equation}
Note that the phase speed of the linear part of \eqref{E:KdV-nonlinear} agrees with $c_0$ in \eqref{E:dispersion2}. 

Therefore, for small amplitude and shallow water waves, satisfying
\begin{equation}\label{def:KdV regime}
a/h_0=(h_0/\ell)^2\ll1,
\end{equation} 
we combine \eqref{E:KdV-linear} and \eqref{E:KdV-nonlinear}, to arrive at {\em the vorticity-modified Korteweg-de Vries equation}:
\begin{equation}\label{E:KdV}
\eta_t+c_0\eta_x+c_2\eta_{xxx}+\frac12\frac{3g+\omega^2h_0}{\sqrt{gh_0+\frac14\omega^2h_0^2}}\eta\eta_x=0,
\end{equation}
where $c_0$ and $c_2$ are in \eqref{E:dispersion2}. In the irrotational setting, \eqref{E:KdV} becomes 
\begin{equation}\label{E:KdV0}
\eta_t+\sqrt{gh_0}(\eta_x+h_0^2\eta_{xxx})+\frac32\sqrt{\frac{g}{h_0}}\eta\eta_x=0,
\end{equation}
and agrees with \cite[(13.99)]{Whitham}, for instance. 

One may modify the argument in \cite[Section~7.1]{Lannes}, for instance, in the irrotational setting, to rigorously justify \eqref{E:KdV}. But we do not include the details~here. 

Constant vorticity does not qualitatively change wave breaking in the Korteweg-de Vries theory. As a matter of fact, for any $\omega\in\mathbb{R}$, the Cauchy problem associated with \eqref{E:KdV} is well-posed globally in time in $H^1(\mathbb{R})$, say, similarly to the irrotational setting.

\subsection*{Extensions}

The preceding argument may readily be adapted to many nonlinear waves in dispersive media, other than surface water waves. Suppose that $c(k)$ represents the phase speed in the linear theory. For long waves satisfying $|k|\ll 1$, where $k$ is the wave number, one may expand $c$ up to terms order $k^2$, provided that it is real valued, smooth, and even, to find 
\[
c(k)\sim c_0-c_2k^2
\]
for some constants $c_0$ and $c_2$. A simplest linear equation whose dispersion relation is $c_0-c_2k^2$ would be 
\[
\eta_t+c_0\eta_x+c_2\eta\eta_{xxx}=0.
\]
It is then only necessary to have access to the form of the nonlinearity. For small amplitude waves satisfying $a\ll1$, where $a$ is the amplitude, one typically finds the nonlinear term in \eqref{E:KdV-nonlinear} or \eqref{E:KdV} at leading order. After normalization of parameters, therefore, we arrive at
\[
\eta_t+c_0\eta_x+c_2\eta_{xxx}+\eta\eta_x=0.
\]
To recapitulate, the Korteweg-de Vries equation models many nonlinear waves in dispersive media in a small amplitude and long wave regime. 
 
\section{The Whitham equation}\label{sec:Whitham}

As Whitham~\cite{Whitham} emphasized, ``the breaking phenomenon is one of the most intriguing long-standing problems of water wave theory." Recall from Section~\ref{sec:shallow} that \eqref{E:shallow} explains wave breaking. But the shallow water theory goes too far. It predicts that all solutions carrying an increase of elevation break. Yet it is a matter of experience that some waves in water do not break. Perhaps, the dispersion effects, which drop out in \eqref{E:shallow}, inhibit breaking.

By the way, when gradients are no longer small, the shallow water assumption \eqref{def:shallow water}, under which one advocates \eqref{E:shallow}, is no longer valid, and \eqref{E:shallow solution}-\eqref{E:char} loses relevance well before breaking occurs. Nevertheless, as Whitham argued, ``breaking certainly does occur and in some circumstances does not seem to be too far away from" what \eqref{E:shallow solution}-\eqref{E:char} describes. 

But recall from the previous section that, including some dispersion effects, \eqref{E:KdV} goes too far and predicts that no solutions break. Therefore, one necessitates some dispersion effects to satisfactorily explain breaking, but the dispersion of the Korteweg-de Vries equation seems too strong. This is not surprising because \eqref{E:dispersion2} poorly approximates \eqref{E:dispersion} when $kh_0$ becomes large.

Whitham~\cite{Whitham} noted that ``it is intriguing to know what kind of simpler mathematical equations (than the governing equations of the water wave problem) could include" the breaking effect, and in the irrotational setting, he put forward
\begin{equation}\label{E:Whitham0}
\eta_t+c(|\partial_x|)\eta_x+\frac32\sqrt{\frac{g}{h_0}}\eta\eta_x=0,
\end{equation}
where $c(|\partial_x|)$ is a Fourier multiplier operator, defined as 
\begin{equation}\label{def:c0}
\widehat{c(|\partial_x|)f}(k)=\sqrt{\frac{g\tanh(|k|h_0)}{|k|}}\widehat{f}(k)
\end{equation}
in a suitable function space. Here and elsewhere, the circumflex means the Fourier transform.
This combines the dispersion relation in the linear theory of water waves (see \eqref{E:dispersion0}) and a nonlinearity of the shallow water theory (see \eqref{E:KdV0}) in the irrotational setting. For small amplitude and shallow water waves satisfying \eqref{def:KdV regime}, the solutions of \eqref{E:Whitham0}-\eqref{def:c0} differ from those of \eqref{E:KdV0} merely by higher order terms during the relevant time scale; see \cite[Section~7.4]{Lannes}, for instance, for details. But the Whitham equation may offer improvements over the Korteweg-de Vries equation for relatively deep water or, equivalently, relatively short waves. As a matter of fact, numerical experiments (see \cite{MKD}, for instance) indicate that the Whitham equation approximates solutions of the water wave problem on par with or better than the KdV equation or other shallow water models in some respects, outside of the long wave regime. Whitham conjectured wave breaking in \eqref{E:Whitham0}-\eqref{def:c0}. The author~\cite{Hur-breaking} recently solved it. 

For constant vorticities, we take matters further and propose 
\begin{equation}\label{E:Whitham}
\eta_t+c(|\partial_x|)\eta_x+\frac12\frac{3g+\omega^2h_0}{\sqrt{gh_0+\frac14\omega^2h_0^2}}\eta\eta_x=0,
\end{equation}
where 
\begin{equation}\label{def:c}
\widehat{c(|\partial_x|;\omega)f}(k)=\Bigg(\frac{\omega}{2}\frac{\tanh(|k|h_0)}{|k|}
+\sqrt{\frac{g\tanh(|k|h_0)}{|k|}+\frac{\omega^2}{4}\frac{\tanh^2(|k|h_0)}{k^2}}\Bigg)\widehat{f}(k).
\end{equation}
This combines the dispersion relation in the linear theory of water waves (see \eqref{E:dispersion}) and a nonlinearity of the shallow water theory (see \eqref{E:KdV}) in the constant vorticity setting. For small amplitude and shallow water waves, satisfying \eqref{def:KdV regime}, one may modify the argument in \cite[Section~7.4]{Lannes}, for instance, in the irrotational setting, to verify that \eqref{E:Whitham}-\eqref{def:c} is equivalent to \eqref{E:KdV} during the relevant time scale. But the Whitham equation may offer improvements over the Korteweg-de Vries equation for relatively deep water or, equivalently, relatively short waves, similarly to the irrotational setting. 

\subsection*{Wave breaking}

One may follow along the same line as the proof in \cite{Hur-blowup}, for instance, to establish a H\"older norm blowup in \eqref{E:Whitham}-\eqref{def:c}. Specifically, if $\eta_0\in L^2(\mathbb{R})\bigcap C^{1+\alpha}(\mathbb{R})$, $0<\alpha<1$, and if for any $0<\delta<1$ and $2<p<3$, 
\[
\int_{|x|<1}(\eta_0(x)-\eta_0(0))\text{sgn}(x)|x|^{-\delta}~dx
+\int_{|x|>1}(\eta_0(x)-\eta_0(0))\text{sgn}(x)|x|^{-p}~dx>C
\]
for some constant $C>0$, then the solutions of \eqref{E:Whitham}-\eqref{def:c} and $\eta(\cdot,0)=\eta_0$ exhibits
\begin{equation}\label{E:H-blowup}
\lim_{t\to t_*}\|\eta(\cdot,t)\|_{C^{1+\alpha}}(\mathbb{R})=\infty
\end{equation}
for some $t_*>0$. 

To see this, we may rewrite \eqref{E:Whitham}-\eqref{def:c}, after normalization of parameters, as
\[
\eta_t+\mathscr{H}|\partial_x|^{1/2}(1+\mathscr{R}(|\partial_x|))\eta+\eta\eta_x=0,
\]
where $\mathscr{H}$ is the Hilbert transform and 
\[
|\mathscr{R}(|k|)|\leq e^{-|k|}\qquad\text{pointwise in $\mathbb{R}$}.
\]
The author~\cite{Hur-blowup} proved the H\"older norm blowup in 
\[
\eta_t+\mathscr{H}|\partial_x|^{1/2}\eta+\eta\eta_x=0.
\]
We then observe that the proof remains valid under perturbations of $\mathscr{H}|\partial_x|^{1/2}$ by linear operators bounded in $L^2(\mathbb{R})$. 

Furthermore, one may promote \eqref{E:H-blowup} to wave breaking. Specifically, if $\eta_0\in H^\infty(\mathbb{R})$,
\[
\epsilon^2(\inf_{x\in\mathbb{R}}\eta_0'(x))^2>1+\|\eta\|_{H^3(\mathbb{R})}
\]
for $\epsilon>0$ sufficiently small, and $\eta_0$ belongs to the Gevrey class of index $2$, then the solution of \eqref{E:Whitham0}-\eqref{def:c0} and $\eta(\cdot, 0)=\eta_0$ exhibits wave breaking. That is,
\[
|\eta(x,t)|<\infty \qquad\text{for all $x\in\mathbb{R}$}\quad\text{for all $t\in[0,t_*)$}
\]
but
\[
\inf_{x\in\mathbb{R}}\eta_x(x,t)\to-\infty\qquad\text{as $t\to t_*-$}
\]
for some $t_*>0$. 

To see this, we follow along the same line as the proof in \cite{Hur-breaking}, but we use 
\[
\|c(|\partial_x|)\partial_x^n\eta_x\|_{L^\infty(\mathbb{R})}\leq 
C\||\partial_x|^n\eta\|_{L^\infty(\mathbb{R})}^{1/2}\||\partial_x|^{n+1}\eta\|_{L^\infty(\mathbb{R})}^{1/2}
\]
by the Gagliardo-Nirenberg interpolation inequality, in place of a quantitative study of the Fourier transform of $c(|\partial_x|)$. Better yet, we may establish wave breaking in equations of Whitham type for a wide range of dispersion relation. 

\section{The full-dispersion shallow water equations}\label{sec:MI}

In the 1960s, Benjamin and Feir~\cite{BF, BH} and Whitham~\cite{Whitham1967} discovered that a Stokes wave would be unstable to long wavelength perturbations --- namely, the Benjamin-Feir or modulational instability --- provided that 
\[
kh_0>1.363\dots.
\]
Corroborating results arrived about the same time, but independently, by Lighthill \cite{Lighthill}, Zakharov \cite{Zakharov-WW}, among others. This is argued to eventually lead to wave breaking; see \cite{BH, LHCokelet1978}, for instance. As a matter of fact, experimental studies (see \cite{Melville1982, Melville1983}, for instance) confirm a ``frequency downshift" when the wave at maximum modulation is breaking or close to breaking. In the 1990s, Bridges and Mielke~\cite{BM1995} addressed the corresponding spectral instability in a rigorous manner. By the way, it is difficult to justify the 1960s theory in a functional analytic setting. But the proof does not easily permit the effects of surface tension or constant vorticity. 

In recent years, the Whitham equation has gathered renewed attention because of its ability to explain high frequency phenomena of water waves; see \cite{Hur-breaking} and references therein. In particular, Johnson and the author~\cite{HJ2} demonstrated that a small amplitude and $2\pi/k$ periodic traveling wave of \eqref{E:Whitham0}-\eqref{def:c0} be spectrally unstable to long wavelength perturbations, provided that $kh_0>1.145\dots$, similarly to the Benjamin-Feir instability, and it is stable to square integrable perturbations otherwise. By the way, \eqref{E:shallow0} does not admit periodic traveling waves, and periodic wavetrain of \eqref{E:KdV0} are all modulationally stable. 

Johnson and the author~\cite{HJ3} took matters further and included the effects of surface tension, by replacing \eqref{def:c0} by 
\begin{equation}\label{def:cT}
\widehat{c(|\partial_x|)f}(k)=\sqrt{(g+Tk^2)\frac{\tanh(|k|h_0)}{|k|}}\widehat{f}(k)
\end{equation}
(see \eqref{E:T-dispersion0}), and the effects of constant vorticity by replacing \eqref{def:c0} by \eqref{def:c}. The results by and large agree with \cite{Kawahara, DR} and \cite{TKM}, for instance, from multiple scale expansions of the physical problem. But including the effects of surface tension, it fails to predict the critical wave number at the ``large surface tension" limit. Including the effects of constant vorticity, it fails to predict the stability at some ``large vorticity" limits. Recently, Pandey and the author~\cite{HP2,HP3} extended the Whitham equation to include higher order nonlinearities and to permit bidirectional propagation, respectively, which correctly predict the capillary effects on the Benjamin-Feir instability.

Including the effects of constant vorticity, here we propose {\em the full-dispersion shallow water equations}:
\begin{equation}\label{E:FDSW}
\begin{aligned}
&\eta_t+(u(h_0+\eta))_x-\omega \eta\eta_x=0,\\
&u_t+\omega h_0c^2(|\partial_x|)u_x+c^2(|\partial_x|)\eta_x+uu_x=0,
\end{aligned}
\end{equation}
where 
\begin{equation}\label{def:c^2}
\widehat{c^2(|\partial_x|)f}(k)=\frac{g\tanh(|k|h_0)}{|k|}\widehat{f}(k).
\end{equation}
They combine the dispersion relation in the linear theory of water waves in the constant vorticity setting and the vorticity-modified shallow water equations, and they extend \eqref{E:Whitham}-\eqref{def:c} to permit bidirectional propagation. As a matter of fact, the phase speed of the linear part of \eqref{E:FDSW}-\eqref{def:c^2} satisfies \eqref{E:dispersion}. In the limit as $kh_0\to0$, \eqref{E:FDSW}-\eqref{def:c^2} becomes \eqref{E:shallow}. Moreover, for small amplitude and long waves, after normalization of parameters, \eqref{E:FDSW}-\eqref{def:c^2} has relevance to a variant of the Boussinesq equations (see \cite{YGT1994}, for instance)
\begin{equation}\label{E:Boussinesq}
\begin{aligned}
&\eta_t+u_x+(u\eta)_x-\omega \eta\eta_x=0,\\
&u_t-\frac13u_{xxt}+\frac13\omega u_{xxx}+\eta_x+uu_x=0,
\end{aligned}
\end{equation}
In the irrotational setting, \eqref{E:FDSW}-\eqref{def:c^2} becomes
\begin{equation}\label{E:FDSW0}
\begin{aligned}
&\eta_t+(u(h_0+\eta))_x=0,\\
&u_t+c^2(|\partial_x|)\eta_x+uu_x=0,
\end{aligned}
\end{equation}
and agrees with what \cite{HP2} proposes. 

\subsection*{Modulational instability}

In what follows, we assume that $g=h_0=1$ for simplicity of notation.

One may follow along the same line as the argument in \cite{HP2} and establish that for any $\omega\in\mathbb{R}$ for any $k>0$, a one parameter family of small amplitude and $2\pi/k$-periodic traveling wave of \eqref{E:FDSW}-\eqref{def:c^2} exists, denoted $\eta(a)(z), u(a)(z)$ and $c(a)$, for $a\in\mathbb{R}$ and $|a|$ sufficiently small for $z=k(x-c(a)t)$. We omit the details. A straightforward calculation reveals that 
\begin{align*}
\eta(a)(z)=&a\cos z+\frac14a^2\Big(
\frac{3c_0^2-3\omega c_0+\omega^2}{c_0^2-\omega c_0-1}\\
&\hspace*{65pt}+\frac{3c_0-\omega c_0(1+2c^2(2k))+\omega^2c^2(2k)}{c_0^2-\omega c^2(2k)c_0-c^2(2k)}\cos(2z)\Big)+O(a^3),\\
u(a)(z)=&ac_0\cos z+\frac14a^2\Big(
\frac{c_0^3+2c_0-\omega}{c_0^2-\omega c_0-1}\\
&\hspace*{70pt}+\frac{c_0^3+2c^2(2k)c_0+\omega c^2(2k)}{c_0^2-\omega c^2(2k)c_0-c^2(2k)}\cos(2z)\Big)+O(a^3),
\end{align*}
and $c(a)=c_0+O(a^2)$, where
\begin{equation}\label{E:c0}
c_0=\frac{\omega}{2}\frac{\tanh(k)}{k}+\sqrt{\frac{\tanh(k)}{k}+\frac{\omega^2}{4}\frac{\tanh^2(k)}{k^2}}
\end{equation}
satisfies \eqref{E:dispersion}.

One may follow along the same line as the argument in \cite{HP2} (see also \cite{HJ2, HJ3, BHJ}) and make a very lengthy and complicated, albeit explicit, spectral perturbation calculation for the linearized operator associated with \eqref{E:FDSW}-\eqref{def:c^2}. The result states that a small amplitude and $2\pi/k$-periodic traveling wave of \eqref{E:FDSW}-\eqref{def:c^2} is unstable to long wavelength perturbations, provided that 
\begin{equation}\label{def:ind}
\text{ind}(k)=\frac{i_1(k)i_2(k)}{i_3(k)}i_4(k)<0,
\end{equation}
where
\begin{align*}
i_1(k)=&(kc(k))'',\\
i_2(k)=&((kc(k))')^2-c^2(0), \\
i_3(k)=&c^2(k)-c^2(2k),
\end{align*}
and $i_4$ is an explicit function of $k$, but it is so lengthy and complicated that we do not include here. Moreover, the small amplitude and periodic traveling wave is spectrally stable to square integrable perturbation in the vicinity of the origin in $\mathbb{C}$ if $\text{ind}(k)>0$. We omit the details.

The modulational instability index reveals four resonance mechanisms which change the sign of \eqref{def:ind} and, hence, the modulational stability and instability for \eqref{E:FDSW}-\eqref{def:c^2}, similarly to the irrotational setting. Note that
\[
c(k)=\text{the phase velocity}\quad\text{and}\quad (kc(k))'=\text{the group velocity}
\]
in the linear theory. Specifically, 
\begin{itemize}
\item[(R1)] $i_1(k)=0$ at some $k$; the group velocity achieves an extremum at the wave number $k$; 
\item[(R2)] $i_2(k)=0$ at some $k$; the group velocity at the wave number $k$ coincides with the phase velocity in the limit as $k\to0$, resulting in the ``resonance of short and long waves";
\item[(R3)] $i_3(k)=0$ at some $k$; the phase velocities of the fundamental mode and the second harmonic coincide at the wave number $k$, resulting in the ``second harmonic resonance";
\item[(R4)] $i_4(k)=0$ at some $k$, resulting in resonance of the dispersion and nonlinear effects.
\end{itemize}

\begin{figure}[h]
\includegraphics[scale=0.7]{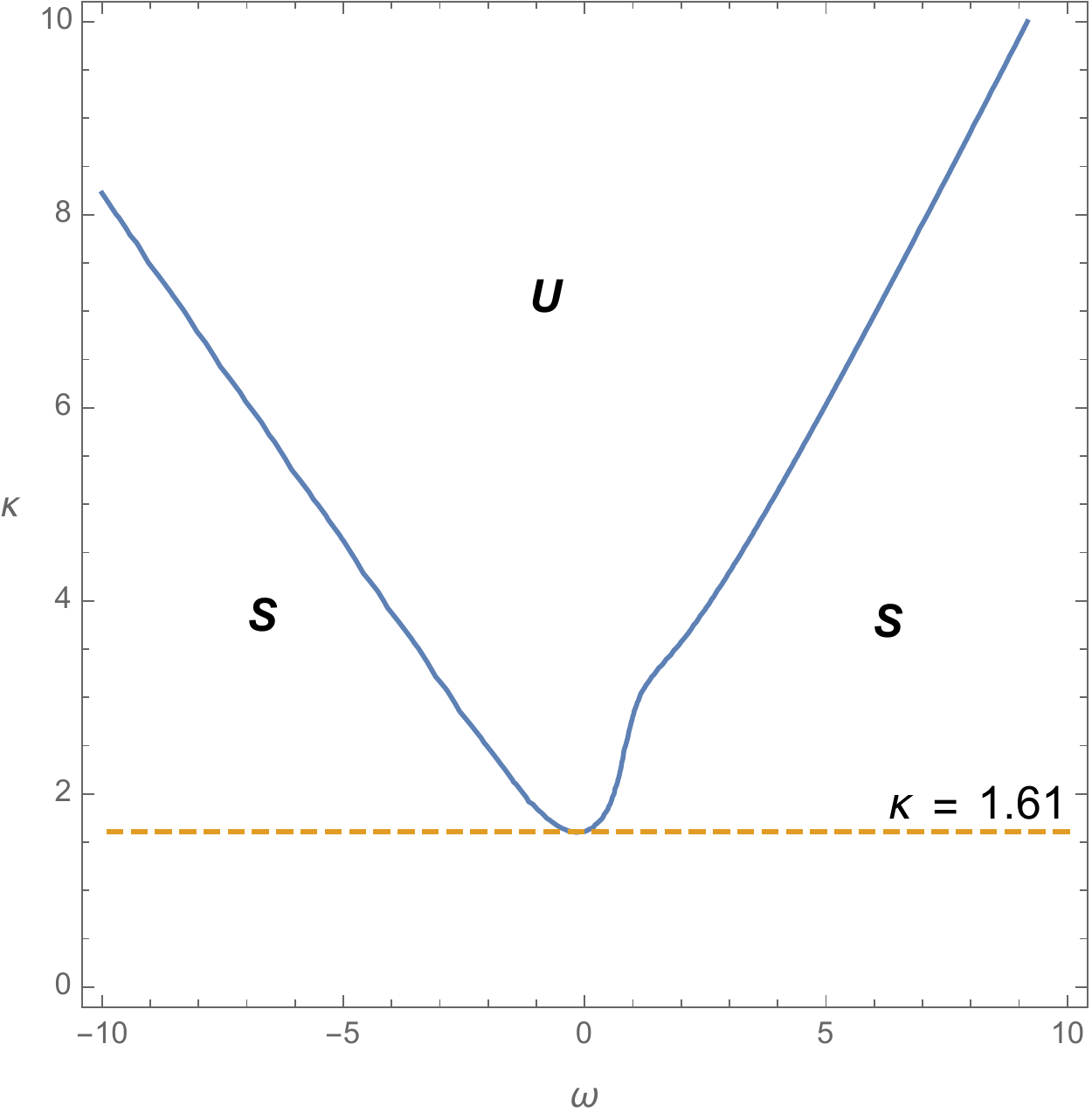}
\caption{Modulational instability diagram in the $\omega$ versus $k$ plane for a small amplitude and periodic traveling wave of \eqref{E:FDSW}-\eqref{def:c^2}. ``S'' and ``U'' denote the regions of modulational stability and instability, respectively. In the case of $\omega=0$, modulational stability and instability changes when $k=1.610\dots$.}\label{fig:T=0}
\end{figure}

For any $\omega\in\mathbb{R}$, a straightforward calculation reveals that resonances (R1) through (R3) do not occur. A numerical evaluation reveals that $i_4$ changes its sign once, from positive to negative, at $k_c(\omega)$, say, depending on $\omega$. To summarize, a small amplitude and $2\pi/k$-periodic traveling wave of \eqref{E:FDSW}-\eqref{def:c^2} is unstable to long wavelength perturbations if $k>k_c(\omega)$, and it is stable if $k<k_c(\omega)$. In the irrotational setting, $k_c(0)=1.610\dots$, and agrees with that in \cite{HP2}. 

Therefore, constant vorticity does not seem to qualitatively change modulational stability and instability. But, a numerical evaluation reveals that $k_c(\omega)$ increases unboundedly as $|\omega|\to\infty$; see Figure~\ref{fig:T=0}. In particular, for any $k>0$, a small amplitude and $2\pi/k$ periodic traveling wave of \eqref{E:FDSW}-\eqref{def:c^2} is modulationally stable, provided that $|\omega|$ is sufficiently large, depending on $k$ and the sign of $\omega$. To compare, the result in \cite{HJ3} based on the Whitham equation fails to predict the stabilizing effects of large vorticities, both positive and negative.

\begin{figure}[h]
\includegraphics[scale=0.7]{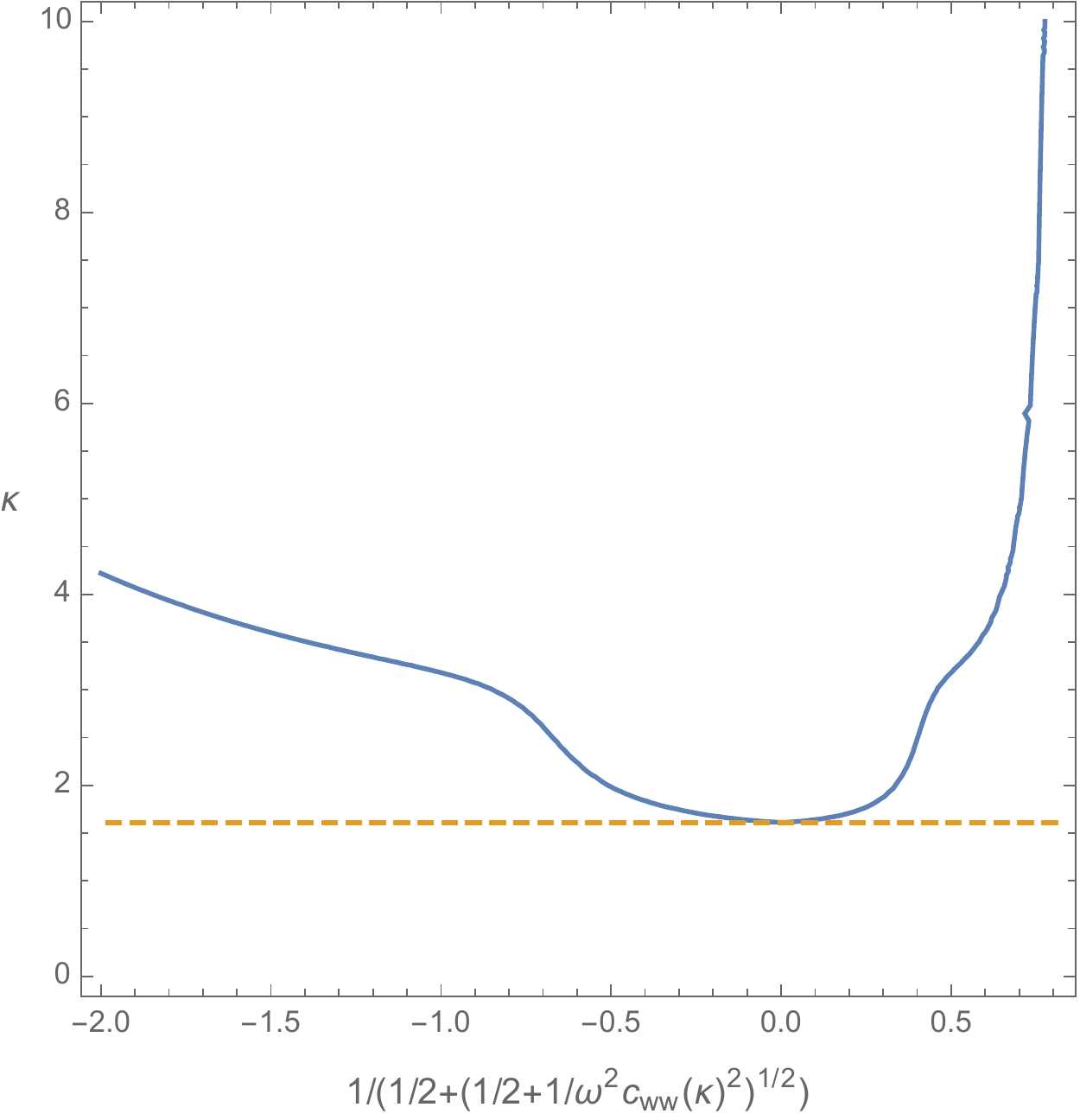}
\caption{Modulational instability diagram in the $\omega(c_0(k)/c(k))$ versus $k$ plane for a small amplitude and periodic traveling wave of \eqref{E:FDSW}-\eqref{def:c^2}. The curve is asymptotic to the vertical at $\omega(c_0(k)/c(k))=0.807\dots$.}\label{fig:T=01}
\end{figure}

Moreover, the result qualitatively agrees with that in \cite{TKM2012}, for instance, from multiple scale expansion of the physical problem; compare Figure~\ref{fig:T=01} and \cite[Figure~3]{TKM2012}. The critical renormalized vorticity $=\omega(c_0(k)/c(k))$, for which small amplitude and periodic traveling waves are all modulationally stable is $0.807\dots$. This compares reasonably well with $2/3$ in \cite{TKM2012}.

\subsection*{Effects of constant vorticity and surface tension}

Including the effects of constant vorticity and surface tension, we replace \eqref{def:c^2} by 
\begin{equation}\label{def:c^2T}
\widehat{c^2(|\partial_x|)f}(k)=(g+Tk^2)\frac{\tanh(|k|h_0)}{|k|}\widehat{f}(k).
\end{equation}

\begin{figure}[h]
\includegraphics[scale=0.7]{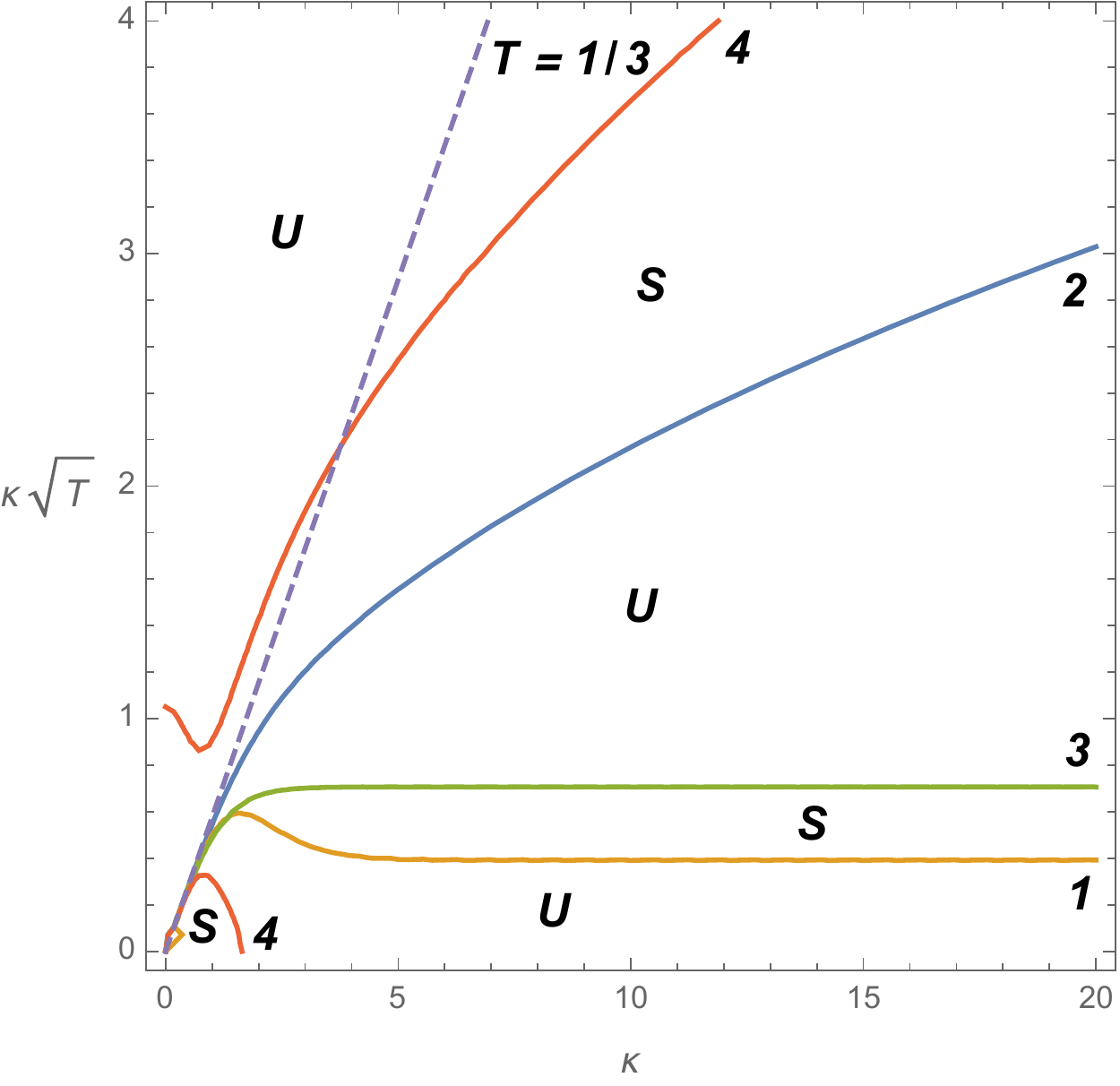}
\caption{Mosulational instability diagram in the $k$ versus $k\sqrt{T}$ plane for small amplitude and periodic traveling wave of \eqref{E:FDSW0} and \eqref{def:c^2T} when $\omega=0$. To interpret, for any $T>0$ one must invasion a line through the origin with the slope $T$. ``S'' and ``U'' denote the regions of modulational stability and instability, respectively. Solid curves represents roots of the modulational instability index and are labeled according to their mechanism.}\label{fig:w=0}
\end{figure}

Of course, surface tension enters nonlinearly in the water wave problem and may contribute to higher order nonlinearities in \eqref{E:FDSW}. As a matter of fact, \eqref{E:Whitham} and \eqref{def:cT} fail to predict the critical wave number at the ``large surface tension" limit (see \cite{HJ3}), whereas a model combining a Camassa-Holm equation and \eqref{def:cT} correctly predicts the limit (see \cite{HP3}). Nevertheless, in the irrotational setting, \eqref{E:FDSW0} and \eqref{def:c^2T} correctly explain the modulational stability and instability of gravity capillary waves (see \cite{HP2}).

Specifically, when $\omega=0$, for $0<T<1/3$, resonances (R1) through (R4) all occur, but at distinct wave numbers, resulting in intervals of modulationally stable and unstable wave numbers. For $T>1/3$, the capillary effects do not qualitative change the modulational stability and instability for \eqref{E:FDSW}. We summarize the result in Figure~\ref{fig:w=0}. The result qualitatively agrees with that in \cite{Kawahara, DR}, for instance, from multiple scale expansion of the physical problem. Compare Figure~\ref{fig:w=0} and \cite[Figure~1]{DR}, for instance.

\begin{figure}[h]
\includegraphics[scale=0.7]{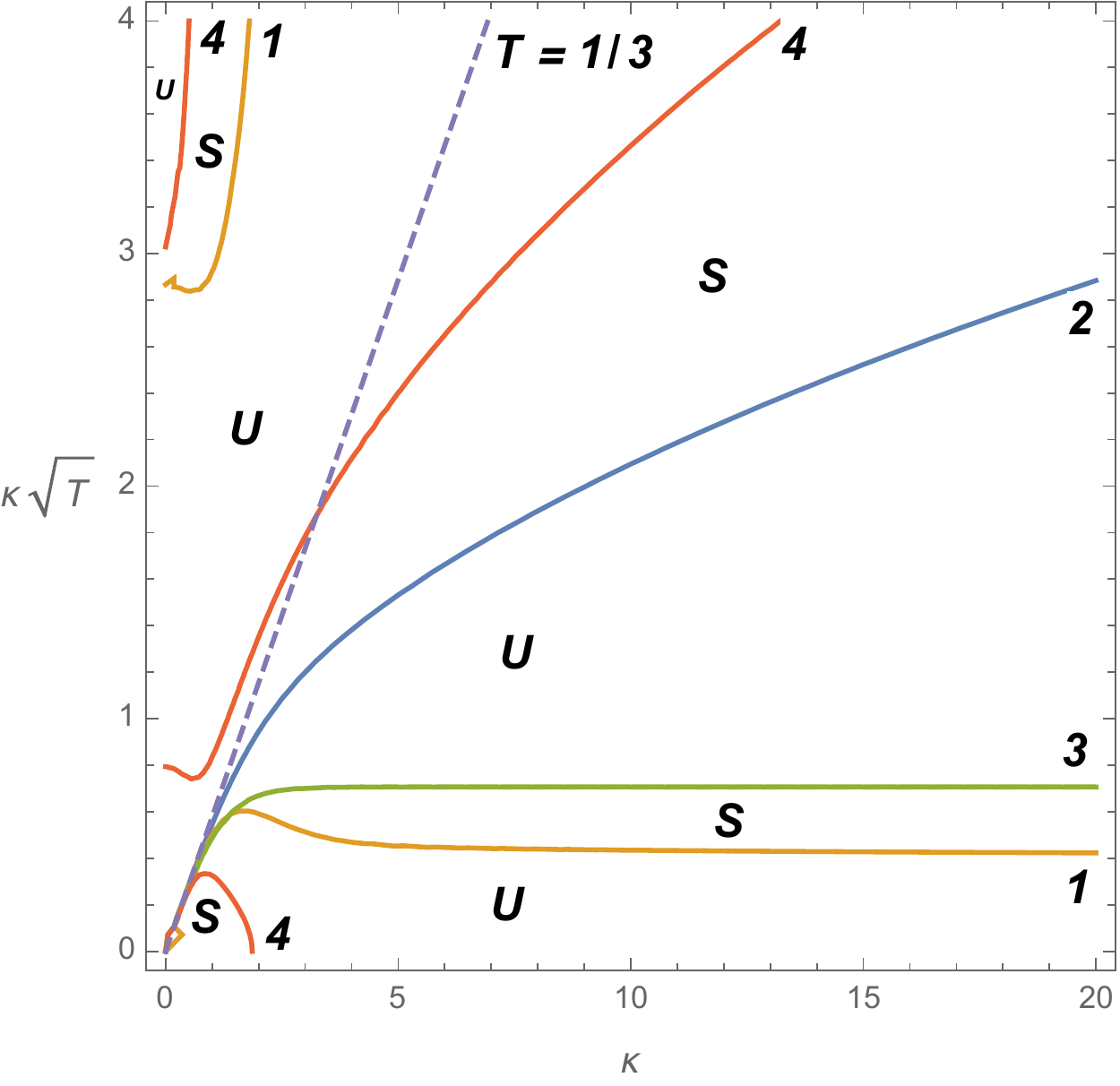}
\caption{Modulational instability diagram in the $k$ versus $k\sqrt{T}$ plane for small amplitude and periodic traveling wave of \eqref{E:FDSW} and \eqref{def:c^2T} when $\omega=-3$. }\label{fig:w=-3}
\end{figure}

\begin{figure}[h]
\includegraphics[scale=0.7]{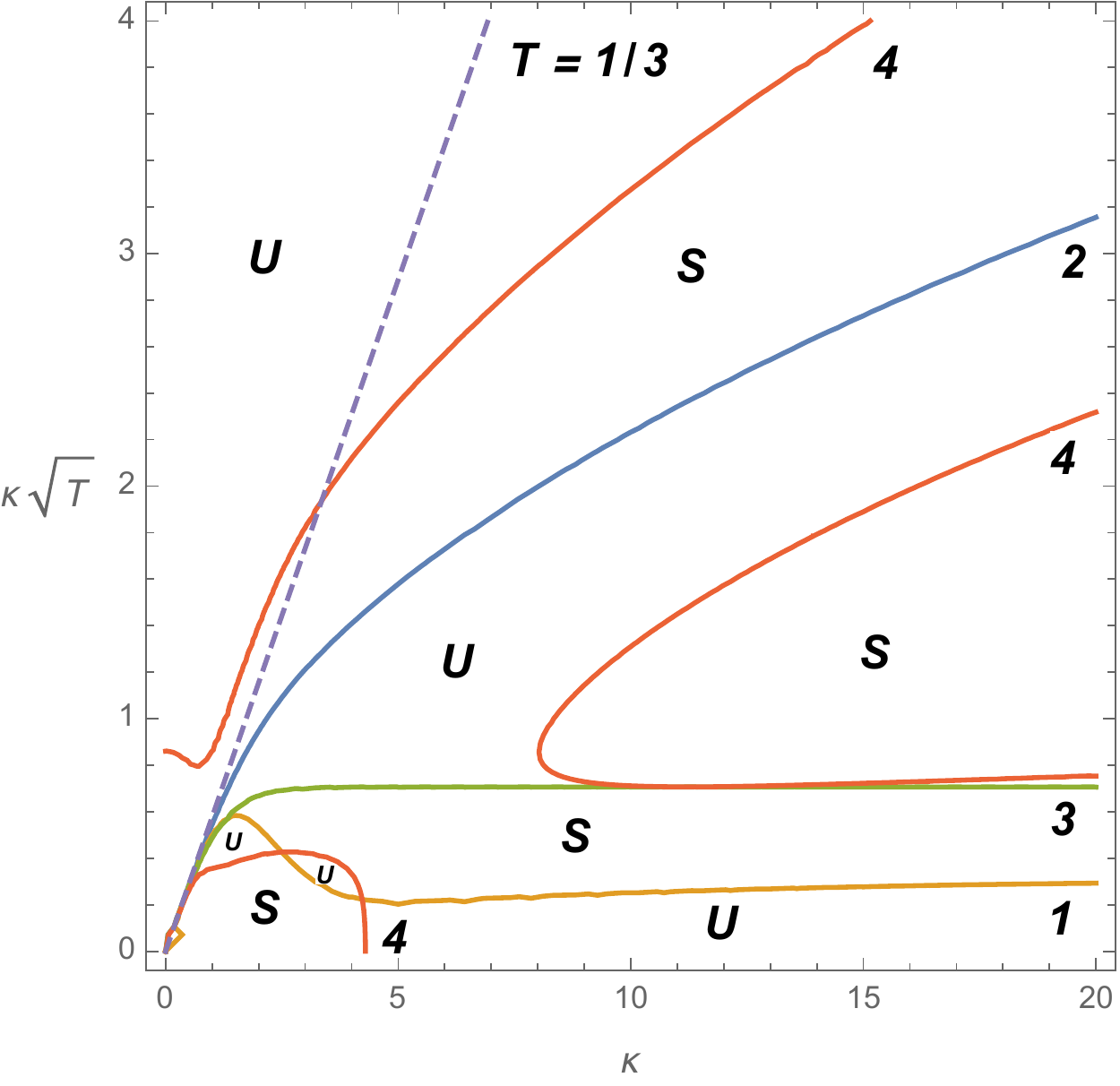}
\caption{Modulational instability diagram in the $k$ versus $k\sqrt{T}$ plane for small amplitude and periodic traveling wave of \eqref{E:FDSW} and \eqref{def:c^2T} when $\omega=3$.}\label{fig:w=3}
\end{figure}

To proceed, for nonzero constant vorticities, we find a very rich set of resonant wave numbers, for which resonances (R1)-(R4) occur. Moreover, for sufficiently positive vorticities, we find some wave numbers, for which a linear resonance (R1) or (R2) and the nonlinear resonance (R4) both occur. Therefore, constant vorticity considerably alters the modulational stability and instability for \eqref{E:FDSW} and \eqref{def:c^2T}; see Figure~\ref{fig:w=-3} and Figure~\ref{fig:w=3}, for instance. It is interesting to compare the result with that from a multiple scale expansion of the physical problem, which extends \cite{TKM}, for instance, in the presence of the effects of surface tension.

\subsection*{Acknowledgements}
The author thanks the organizers of the workshop ``B'Wave 2016: A workshop focussed on wave breaking in oceanic and coastal waters" and University of Bergen in Norway. She thanks the anonymous referees for helpful comments and suggestions.

The author is supported by the National Science Foundation under the Faculty Early Career Development (CAREER) Award DMS-1352597, an Alfred P. Sloan Research Fellowship, a Simons Fellowship in Mathematics, and by the University of Illinois at Urbana-Champaign under the Arnold O. Beckman Research Awards Nos. RB14100 and RB16227. She is grateful to the Department of Mathematics at Brown University and  for its generous hospitality. 

\bibliographystyle{amsalpha}
\bibliography{BW}

\end{document}